\documentclass[11pt]{amsart}
\usepackage{amsmath,amssymb,latexsym,dsfont}
\usepackage[left=2.6cm,right=2.6cm,top=2.7cm,bottom=2.7cm]{geometry}

\usepackage{amsmath,amsfonts,amssymb,epsfig, textcomp}
\usepackage{graphicx,tikz}
\usepackage{hyperref, cleveref} 
\usepackage{cases,listings}
\usepackage{color}
\usepackage{mathabx}

\usepackage[nobysame]{amsrefs}

\newtheorem{theorem}{Theorem}[section]
\newtheorem{lemma}{Lemma}[section]

\newtheorem{proposition}{Proposition}[section]

\newtheorem{remark}{Remark}[section]

\newtheorem{definition}{Definition}[section]
\numberwithin{equation}{section}

      \newcommand{\hy}{\hat y}

      \newcommand{\hh}{\hat h}

   \newcommand{\tg}{\widetilde g}
   \newcommand{\ty}{\widetilde y}

      \newcommand{\cN}{{\mathcal N}}

      \newcommand{\N}{\mathbb{N}}

      \newcommand{\eps}{\varepsilon}
      \newcommand{\mR}{\mathbb{R}}

      \newcommand{\mC}{\mathbb{C}}

      \makeatletter
      \def\@setcopyright{}
      \def\serieslogo@{}
      \makeatother

\newcommand{\cL}{\mathcal L}

\newcommand{\be}{\begin{equation}}
\newcommand{\ee}{\end{equation}}

\newcommand{\tf}{{\widetilde f}}

\newcommand{\cD}{{\mathcal D}}

\newcommand{\mH}{\mathbb{H}}

\newcommand{\tM}{\widetilde{M}}

\title[Boundary stabilization of a KdV system]{Rapid and finite-time boundary stabilization of a KdV system}

\author[H.-M. Nguyen]{Hoai-Minh Nguyen}
\address[H.-M. Nguyen]{Laboratoire Jacques Louis Lions, \newline\indent
Sorbonne Universit\'e\newline\indent
Paris, France}
\email{hoai-minh.nguyen@sorbonne-universite.fr}

\begin{document}

\maketitle

\begin{abstract} We construct a static feedback control in a trajectory sense and a dynamic feedback control
to obtain the local rapid boundary stabilization of a KdV system using Gramian operators.  We also construct a time-varying feedback control in the trajectory sense and a time varying dynamic feedback control to reach the local finite-time boundary stabilization for the same system.  
\end{abstract}

\tableofcontents

\noindent {\bf Keywords}: stabilization, rapid stabilization, feedback, dynamic feedback, KdV equation, Grammian operators. 

\noindent{\bf MSC}: 93B52; 93D15; 35B30; 35B35.


\section{Introduction}

This paper is devoted to studying the local rapid boundary stabilization and the local finite-time boundary stabilization of a KdV system. More precisely, 
we investigate the stabilization of the following control system, for $L > 0$,  
\be \label{sys-KdV}
\left\{\begin{array}{cl}
y_t + y_x + y_{xxx} + y y_x = 0 & \mbox{ in } (0, + \infty) \times (0, L), \\[6pt]
y(\cdot, L) =  y (\cdot, 0) = 0, \ \, y_{x} (\cdot, L) - y_x(\cdot, 0) = u & \mbox{ in } (0, + \infty), \\[6pt]
y(0, \cdot) = y_0 (\cdot) & \mbox{ in } (0, L),   
\end{array} \right. 
\ee
where $y_0 \in L^2(0, L)$ is the initial state, $u \in L^2_{loc}[0, + \infty)$ is a control, and $y(t, \cdot) \in L^2(0, L)$ is the state at time $t$.

The controllability of the {\it linearized} system of \eqref{sys-KdV} depends strongly on $L$. It is known from the work of 
Cerpa and Cr\'epeau \cite{CC09-DCDS}
that there is a discrete set of lengths $\cN$ for which the linearized system is not exactly controllable if $L \in \cN$ and the linearized system is exactly controllable otherwise. 
Concerning system \eqref{sys-KdV}, the rapid stabilization of its {\it linearized} system for non-critical lengths has been obtained by Cerpa and Cr\'epeau \cite{CC09-DCDS}. They used the Gramian operators and the analysis involves the optimal control theory as an application of the result of Urquiza \cite{Urquiza05} (see also \cite{Komornik97}). The feedback is thus understood in a weak sense, see, e.g., \cite{WR00, TWX20,Ng-Riccati}. In a very related setting where one controls the Neumann on the right, i.e., one controls $y_x(\cdot, L)$ instead of $y_x(\cdot, L) - y_x(\cdot, 0)$, the local rapid stabilization of the nonlinear system for non-critical lengths was obtained by Coron and L\"u \cite{CL14} using a technique related to the backstepping method. To our knowledge, the extension to the {\it nonlinear} setting using Gramian operators is open, and the local stabilization in finite time is previously out of reach. The goal of this work is to give an answer to this problem. More precisely, dealing with non-critical lengths, we construct a feedback control in 
a trajectory sense, a notion introduced in \cite{Ng-Riccati},  and a dynamic feedback control 
to obtain the local rapid stabilization of \eqref{sys-KdV} using Gramian operators (see \Cref{thm-W} and \Cref{thm-D}). We also construct a time-varying feedback control in the trajectory sense and a time-varying dynamic feedback control to obtain the local finite-time stabilization of \eqref{sys-KdV}.  The ideas in the study of the rapid stabilization are to modify the approach proposed in \cite{Ng-Riccati} to deal with the non-linear term, which cannot be handled by directly using the proposal given there. Concerning the finite-time stabilization, we additionally combine the ideas in \cite{Ng-Riccati} with the ones proposed by Coron and the author in \cite{CoronNg17} in the spirit of \cite{Ng-Schrodinger}. This thus involves the control cost of the linearized system in small time.

\subsection{Statement of the main results on the rapid stabilization}
Define  $A: \cD(A) \subset L^2(0, L) \to L^2(0, L)$ as follows
$$
\cD(A)= \Big\{w \in H^3(0, L); w(L) = w (0) =  0, w_x(L) = w_x (0) \Big\}, 
$$
and 
\be \label{def-A-KdV}
Aw = - w''' - w \mbox{ for } w \in \cD(A). 
\ee
One can check that $A$ is densely defined and closed in the Hilbert space $L^2(0, L)$ equipped with the standard scalar product. Moreover, 
\be \label{A-skew-adjoint}
\mbox{$A$ is skew-adjoint, i.e.,  $\cD(A) = \cD(A^*)$ and $A^* = - A$}, 
\ee
where $A^*$ denotes the adjoint of $A$.  

Let $B: \mR \to \cD(A^*)'$, where $\cD(A^*)'$ is the dual space of $\cD(A^*)$,  be defined by 
\be
\langle B u, w \rangle_{\cD(A^*)', \cD(A^*)} = u w_x(L).  
\ee
Then $B^*:  \cD(A^*) = \cD(A) \to \mR$ is given as follows, for $w \in \cD(A)$,  
\be \label{def-B*-KdV}
B^*w = w_x(L). 
\ee
Then the linearized system of \eqref{sys-KdV} (around the zero state) can be written under the form 
\be \label{sys-Semigroup}
\left\{\begin{array}{lc}
y' = A y + Bu &  \mbox{ in } (0, + \infty), \\[6pt]
y(0) = y_0  
\end{array} \right. 
\ee
(see, e.g., \cite[Section 3]{Ng-Riccati} for the meaning of \eqref{sys-Semigroup}).
One can check that 
\be \label{admissibility}
\mbox{the control operator $B$ is an admissible control operator for all $L > 0$ and $T>0$}, 
\ee 
i.e., for some positive constant $C = C(T, L)$, 
\be \label{admissibility1}
\int_0^{T} |B^*e^{s A^*} z|^2 \, ds \le C \int_0^L |z(x)|^2 \, dx \mbox{ for all } z \in L^2(0, L), 
\ee 
since (see, e.g., \Cref{lem-KdV1})
\be  \label{admissibility1-***}
\|\xi_x(L, \cdot) \|_{L^2(0, T)} \le C \| \xi_0 \|_{L^2(0, L)} \mbox{ for all } \xi_0  \in L^2(0, L),
\ee 
where $\xi \in X_T$ is the unique solution of the system 
\be \label{sys-z-KdV}
\left\{\begin{array}{cl}
\xi_t + \xi_x + \xi_{xxx}  = 0 & \mbox{ in } (0, T) \times (0, L), \\[6pt]
\xi(\cdot, L) =  \xi (\cdot, 0) = 0, \ \, \xi_x (\cdot, L) - \xi_x(\cdot, 0) = 0 &  \mbox{ in } (0, T),  \\[6pt]
\xi(0, \cdot) = \xi_0 (\cdot) & \mbox{ in } (0, L).  
\end{array} \right. 
\ee
Here and in what follows, we denote
\be
X_T = C([0, T]; L^2(0, L)) \cap L^2( (0, T); H^1(0, L)) \mbox{ for } T>0, 
\ee
and 
\be
X_\infty = C([0, +\infty); L^2(0, L)) \cap L^2_{loc}( [0, + \infty); H^1(0, L)). 
\ee
As usual, we denote $\big(e^{t A^*}\big)_{t \ge 0}$ the semigroup generated by $A^*$.  

The controllability of the linearized system of \eqref{sys-KdV} depends strongly on $L$. Denote 
\be \label{def-cN} 
\cN= \left\{L = 2 \pi \sqrt{\frac{k^2 + k l + l^2}{3}}; k, l \in \N_+ \right\}. 
\ee
It is known that the linearized system of \eqref{sys-KdV} given by  
\be \label{sys-LKdV}
\left\{\begin{array}{cl}
y_t + y_x + y_{xxx}  = 0 & \mbox{ in } (0, T) \times (0, L), \\[6pt]
y(\cdot, L) =  y (\cdot, 0) = 0, \ \, y_{x} (\cdot, L) - y_x(\cdot, 0) = u (\cdot) & \mbox{ in } (0, T), \\[6pt] 
y(0, \cdot) = y_0 (\cdot) & \mbox{ in } (0, L),  
\end{array} \right. 
\ee
is exactly controllable for all  (or for some) $T>0$ if and only if $L \not \in \cN$, a result due to Cerpa and Cr\'epeau \cite{CC09-DCDS}. This is equivalent to the fact that for all (or for some) $T>0$, it holds, for some positive constant $C$, 
\be 
\int_0^T |\xi_x(t, L)|^2 \, dt  \ge C \| z_0\|_{L^2(0, L)},
\ee
for all $\xi_0 \in L^2(0, L)$ where $\xi \in X_T$ is the unique solution of system \eqref{sys-z-KdV} if and only if $L \not \in \cN$.  A very closely related work was previously obtained by Rosier \cite{Rosier97}. 

\medskip 
We are ready to introduce the Gramian operators used in our feedback controls. Given $\lambda > 0$, define $Q = Q(\lambda): L^2(0, L) \to L^2(0, L)$ as follows 
\be \label{def-Q-KdV}
\langle Q z_1, z_2 \rangle_{L^2(0, L)} = \int_{0}^{\infty} e^{-2 \lambda s} \langle B^* e^{-s A^*} z_1, B^* e^{-s A^*} z_2 \rangle_{\mR} \, ds \mbox{ for } z_1, z_2 \in L^2(0, L).   
\ee
Here and in what follows, given a Hilbert space $\mH$, we denote $\langle \cdot, \cdot \rangle_{\mH}$ its scalar product and $\cL(\mH)$ the space of all continuous linear applications from $\mH$ to $\mH$ equipped with the standard norm, which is denoted by $\| \cdot \|_{\cL(\mH)}$. 

It is clear that $Q$ is symmetric. It is worth noting that $Q$ is invertible if $L \not \in \cN$ since the linearized system is exactly controllable in small time. An equivalent way to define $Q$ is given by  
\be \label{def-Q-KdV-PDE}
\langle Q z_1, z_2 \rangle_{L^2(0, L)} = \int_{0}^{\infty} e^{-2 \lambda s}  \xi_{1, x} (s, L) \xi_{2, x} (s, L) \, ds \mbox{ for } z_1, z_2 \in L^2(0, L),   
\ee
where $\xi_j \in X_{\infty}$ (with $j=1, \, 2$) is the unique solution of the system 
\be
\left\{\begin{array}{cl}
\xi_{j, t} + \xi_{j, x} + \xi_{j, xxx}  = 0 & \mbox{ in } (0, + 
\infty) \times (0, L), \\[6pt]
\xi_j(\cdot, L) =  \xi_j (\cdot, 0) = 0, \ \, \xi_{j, x} (\cdot, L) - \xi_{j, x}(\cdot, 0) = 0 & \mbox{ in } (0, + \infty),  \\[6pt]
\xi_{j}(0, \cdot) = z_j &  \mbox{ in } (0, L).  
\end{array} \right. 
\ee
One can check, see, e.g., \cite{Ng-Riccati}, that $Q$ satisfies the following {\it important} property:  
\be \label{identity-Op}
A Q + Q A^* - B B^* + 2 \lambda Q = 0,   
\ee 
where \eqref{identity-Op} is understood in the following sense
\begin{multline} \label{identity-Op-meaning} 
\langle Q z_1, A^*z_2 \rangle_{L^2(0, L)} + \langle A^*z_1, Q z_2 \rangle_{L^2(0, L)}  \\[6pt] 
- \langle B^*z_1,  B^*z_2 \rangle_{\mR}   + 2 \lambda \langle Q z_1, z_2 \rangle_{L^2(0, L)}= 0 \quad  \forall \, z_1, z_2 \in \cD(A^*). 
\end{multline}

We are ready to state the rapid stabilization of \eqref{sys-KdV} in the trajectory sense. 

\begin{theorem} \label{thm-W} Assume that $L \not \in \cN$, and let $\lambda > 0$ and $T_0 > 0$. Define $Q = Q(\lambda)$ by \eqref{def-Q-KdV}.  There exists $\eps > 0$ such that for $y_0 \in L^2(0, L)$ with $\| y_0\|_{L^2(0, L)} \le \eps$, there exists a unique weak solution $(y, \ty) \in X_\infty \times X_\infty$ of the system 
\be \label{thm-W-sys}
\left\{\begin{array}{cl}
y_t  + y_x + y_{xxx} + y y_x = 0 &\mbox{ in } (0, + \infty) \times (0, L), \\[6pt]
\ty_t + \ty_x + \ty_{xxx} + 2 \lambda \ty +   \ty y_x = 0 & \mbox{ in } (0, + \infty) \times (0, L), \\[6pt]
y(\cdot, 0) = y (\cdot, L) = 0, y_x(\cdot, L) - y_x(\cdot, 0) = - \ty_x (\cdot, L) & \mbox{ in } (0, + \infty), \\[6pt]
\ty(\cdot, 0) = \ty (\cdot, L) = 0, \ty_x(\cdot, L) - \ty_x(\cdot, 0) = 0 & \mbox{ in } (0, + \infty), \\[6pt]
y(0) = y_0, \quad  \ty(0) = \ty_0: = Q^{-1} y_0 & \mbox{ in } (0, L).  
\end{array} 
\right. 
\ee
Moreover, we have  
\be \label{thm-W-cl1}
\ty (t, \cdot) = Q^{-1}y(t, \cdot) \mbox{ for } t \ge 0, 
\ee
\be \label{thm-W-cl2}
\|y(t, \cdot) \|_{L^2(0, L)}  \le 2 e^{- 2 \lambda t} \| y_0\|_{L^2(0, L)} 
\mbox{ for } t \in [0, T_0],    
\ee
and 
\be \label{thm-W-cl3}
\|y\|_{X_{T_0}}  \le C \| y_0\|_{L^2(0, L)}   
\ee
for some positive constant $C$ depending only on $L$ and $T_0$.  As a consequence of \eqref{thm-W-cl2}, for every $0 < \gamma < \lambda$, there exists $\eps > 0$ such that for every $y_0 \in L^2(0, L)$ with $\| y_0\|_{L^2(0, L)} \le \eps$, it holds 
\be \label{thm-W-cl2-1}
\|y(t, \cdot) \|_{L^2(0, L)}  \le 2 e^{- 2 \gamma t} \| y_0\|_{L^2(0, L)} 
\mbox{ for } t \in [0, + \infty).     
\ee
\end{theorem}

The proof of \Cref{thm-W} is given in \Cref{sect-W-feedback}.

\begin{remark} \rm The definition of the weak solutions in \Cref{thm-W} is given in \Cref{def-KdV2} in \Cref{sect-Pre}. 
\end{remark}

\begin{remark} \label{rem-W} \rm Some comments on \Cref{thm-W} are in orders.  Since $\ty (t, \cdot) = Q^{-1}y(t, \cdot)$ for $t \ge 0$ by \eqref{thm-W-cl1}, the feedback of \eqref{sys-KdV} can be viewed as 
$$
-  \Big( Q^{-1}y(t, \cdot) \Big)_x. 
$$
We only consider this feedback as a static one in a weak sense, {\it the trajectory sense as used in \cite{Ng-Riccati}},  since for $y \in L^2(0, L)$, it is not clear how to give the sense to the action  $-  \Big( Q^{-1}y(t, \cdot) \Big)_x$. Note that our feedback controls given by  $- \ty_x(t, \cdot)$ via  \eqref{thm-W-cl1} are well-defined in the sense of \Cref{thm-W} for all initial data $y_0 \in L^2(0, L)$. This is different from the one given by the optimal control theory used in \cite{Komornik97,Urquiza05}.  See \cite{Ng-Riccati} for more comments on this aspect. 
\end{remark}

We next deal with the dynamic feedback control. In this direction, we prove the following result. 

\begin{theorem}\label{thm-D} Let $L \not \in \cN$, $\lambda > \lambda_0 > 0$,  $\lambda_1  > 0$, $c_0 > 0$, and $T_0 > 0$, and let $Q = Q(\lambda)$ be defined by \eqref{def-Q-KdV}. Assume that 
\be \label{thm-D-lambda}
\lambda_1 - (2 + c_0) \lambda  > 0.       
\ee
There exists $\eps > 0$ such that for $y_0, \ty_0  \in L^2(0, L)$ with $\| y_0\|_{L^2(0, L)}, \| \ty_0\|_{L^2(0, L)} \le \eps$, there exists a unique weak solution $(y, \ty) \in X_\infty \times X_\infty$ of the system 
\be \label{thm-D-sys}
\left\{\begin{array}{cl}
y_t  + y_x + y_{xxx} + y y_x = 0 &\mbox{ in } (0, + \infty) \times (0, L), \\[6pt]
\ty_t + \ty_x + \ty_{xxx} + 2 \lambda \ty - 
\lambda_1 Q^{-1}(y - Q \ty) + \ty y_x  = 0 & \mbox{ in } (0, + \infty) \times (0, L), \\[6pt]
y(\cdot, 0) = y (\cdot, L) = 0, y_x(\cdot, L) - y_x(\cdot, 0) = - \ty_x (\cdot, L) & \mbox{ in } (0, + \infty), \\[6pt]
\ty(\cdot, 0) = \ty (\cdot, L) = 0, \ty_x(\cdot, L) - \ty_x(\cdot, 0) = 0 & \mbox{ in } (0, + \infty), \\[6pt]
y(0) = y_0, \quad  \ty(0) = \ty_0 & \mbox{ in } (0, L).  
\end{array} 
\right. 
\ee
Moreover, 
\be \label{thm-D-cl1}
\| y\|_{X_{T_0}} + \| \ty\|_{X_{T_0}} \le C \lambda \| Q^{-1} \|  \big( \| y(0, \cdot) \|_{L^2(0, L)} + \|\ty(0, \cdot)\|_{L^2(0, L)} \big)
\ee
and
\begin{multline} \label{thm-D-cl2}
\| y(t, \cdot)\|_{L^2(0, L)} + \| \ty(t, \cdot)\|_{L^2(0, L)} \\[6pt] 
\le C \| Q^{-1}\| e^{ - 2 \lambda t}  \big( \| y(0, \cdot) \|_{L^2(0, L)} + \|\ty(0, \cdot)\|_{L^2(0, L)} \big) \mbox{ for } t \in [0, T_0], 
\end{multline}
where $C$ is a positive constant independent of $\lambda$, $\lambda_1$,  $t$, and $(y_0, \ty_0)$. As a consequence of \eqref{thm-D-cl2}, for every $0 < \gamma < \lambda$, there exist $\eps > 0$ and  $C>0$ such that for $y_0, \ty_0  \in L^2(0, L)$ with $\| y_0\|_{L^2(0, L)}$, $\| \ty_0\|_{L^2(0, L)} \le \eps$, it holds 
\be \label{thm-D-cl2-1}
\| y(t, \cdot)\|_{L^2(0, L)} + \| \ty(t, \cdot)\|_{L^2(0, L)} \le C  e^{ - 2 \gamma t}  \big( \| y(0, \cdot) \|_{L^2(0, L)} + \|\ty(0, \cdot)\|_{L^2(0, L)} \big) 
\mbox{ for } t \in [0, + \infty).     
\ee
\end{theorem}

Here and in what follows, we denote $\|Q^{-1} \|_{\cL(L^2(0, L))}$ and $\|Q \|_{\cL(L^2(0, L))}$ by $\| Q^{-1}\|$ and $\| Q\|$ for notational ease. 

\medskip 
The proof of \Cref{thm-D} is given in \Cref{sect-D-feedback}.

\begin{remark} \rm The definition of the weak solutions in \Cref{thm-D} is given in \Cref{def-KdV2} in \Cref{sect-Pre}. 
\end{remark}

\begin{remark} \label{rem-KdV-1} \rm System \eqref{thm-W-sys} is slightly different from the suggestions in \cite{Ng-Riccati}.   If one closely follows the suggestion in \cite{Ng-Riccati}, the equation of $\ty$ in \eqref{thm-W-sys} would be 
\be\label{rem-KdV1-eq}
\ty' + \ty_x + \ty_{xxx} + 2 \lambda \ty +  \ty (Q \ty)_x = 0 \mbox{ in } (0, + \infty) \times (0, L).
\ee
This requires us to make sense of the term $(Q \ty)_x$, which is not clear since $Q$ is only a continuous, linear map from $L^2(0, L)$ into $L^2(0, L)$. We bypass this issue by anticipating the conclusion and replacing $Q \ty$ by $y$ in \eqref{rem-KdV1-eq}. The term $\ty (Q \ty)_x$ becomes $\ty y_x$ as given in \eqref{thm-W-sys}. Similarly, System \eqref{thm-W-sys} is also slightly different from the suggestions in \cite{Ng-Riccati} so that the nonlinear term can be handled. 
\end{remark}

\begin{remark} \rm  In \Cref{thm-W} and \Cref{thm-D}, a new variable $\ty$ is added. Adding a new variable is very natural and has been used a long time ago in the control theory even in finite dimensions for linear control systems, see, e.g.,  \cite[Section 11.3]{Coron07} and \cite[Chapter 7]{Sontag98}.  Coron and Pradly \cite{CP91} showed that there exists a nonlinear system in finite dimensions for which the system cannot be stabilized by static feedback controls but can be stabilized by dynamic feedback ones.  
Dynamic feedback controls of finite dimensional nature, i.e., the complement system is a system of differential equations, have been previously implemented in the infinite dimensions, see, e.g., \cite{CVKB13,Coron-Ng22}. Our new variables are of infinite dimension nature. These are previously proposed in \cite{Ng-Riccati,Ng-Schrodinger} and 
are inspired by the optimal control theory, see e.g. \cite{Flandoli87-LQR,LT91,WR00,Zwart96} and the references therein. 
\end{remark}

We also construct a static feedback control in the trajectory sense and a dynamic feedback control to obtain the local finite-time stabilization of the KdV. These results are given in \Cref{pro-FT} in \Cref{sect-W-feedback} and \Cref{pro-FTD} in \Cref{sect-D-feedback}, respectively. 

\subsection{Related works} In this section, we briefly discuss the local boundary controllability and the local stabilization of the KdV equation. We first deal with the  controllability.  When the controls are $y(\cdot, 0)$, $y(\cdot, L)$,   $y_x(\cdot, L)$, Russell and Zhang \cite{RZ96} proved that the KdV equation is small time,  locally, exactly controllable. The case of the left Dirichlet boundary control ($y(\cdot, L) = y_x(\cdot, L) = 0$ and $y(\cdot, 0)$ is controlled) was investigated by Rosier \cite{Rosier04} (see also \cite{GG08}).  We next discuss the case where one controls the right Neumann boundary, i.e.,  $y(\cdot, 0) = y(\cdot, L) = 0$ and $y_x(\cdot, L)$ is a control. Rosier~\cite{Rosier97} proved  that the KdV system is small time, locally, exactly controllable provided that the length $L$ is not critical,
i.e., $L \notin \cN$, where $\cN$ is also given by \eqref{def-cN}.  
To tackle the control problem for a critical length $L \in \cN$, Coron and Cr\'epeau
introduced the power series expansion method \cite{CC04}. The idea is to take into account the effect of the nonlinear term $y y_x$  absent in the corresponding linearized system. Using this 
method, Coron and Cr\'epeau showed \cite{CC04} (see also \cite[section 8.2]{Coron07}) that the KdV system is small time,  locally,  exactly controllable when the unreachable space of the linearized system is of dimension 1. Cerpa \cite{Cerpa07} and Cr\'epeau and Cerpa \cite{CC09}
developed the analysis in \cite{CC04} to prove that the KdV system  is {\it finite time}, locally, exactly 
controllable for other critical lengths. With Coron and Koenig \cite{CKN-24}, we proved that such a system is not small time, locally, null controllable for a class of critical lengths.  This fact is surprising when compared with known results on internal controls for the KdV equation. It is known, see \cite{CPR15, PVZ02, Pazoto05}, that the KdV system with $y(\cdot, 0) = y(\cdot, L) = y_x(\cdot, L) = 0$ is small time,  locally controllable using internal controls {\it whenever} the control region contains an {\it arbitrary}, open subset of $(0, L)$.  It is worth noting that without controls, i.e. the control is taken to be zero, the decay of the solutions for critical lengths might occur but very slow, see e.g., \cite{CCS15,TCSC18,Ng-Decay21}. A related control setting is the one where one controls the Dirichlet on the right. This control problem was first investigated by Glass and Guerrero \cite{GG10}. To this end, in the spirit of Rosier's work mentioned above, they introduced the corresponding set of critical lengths, which is some how more involved. 
Concerning such a system, Glass and Guerrero proved that the corresponding linearized KdV system is small time, exactly controllable if $L \not \in \cN_D$. Developing this result, they also established that the KdV system \eqref{sys-KdV}  is small-time locally controllable. Recently, the critical case was handled in \cite{Ng-KdV-D}. To this end, we showed that the KdV system with the Dirichlet controls on the right is not locally null controllable in small time and established that the unreachable space of the linearized system is always of dimension 1. We also provide a criterion for the local controllability in finite time. In particular, we show that there are critical lengths for which the system is not locally null controllable in small time but locally exactly controllable in finite time. These phenomena are quite distinct in comparison with the setting where one controls the Neumann on the right mentioned above. 

\medskip 

The stabilization of the KdV equation has been previously studied with internal controls in \cite{RZ96,MVZ-02,RZ06,MMP07,LRZ10, MCPA17} and the references therein. Concerning the boundary controls for the KdV equation, in addition to the work \cite{CC13, CL14} mentioned previously,  we refer \cite{Zhang94,CC09-DCDS,CRX17} and the references therein. It is worth noting that the backstepping related technique used in \cite{CL14} has been developed to study the stabilization for other settings such as hyperbolic systems \cite{CVKB13,BC16}, wave equations \cite{KGBS08,SCK10},  heat equations \cite{Liu03,CoronNg17}, 
Kuramoto-Sivashinsky equations \cite{CL15}, water waves systems \cite{CHXZ22}, Gribov operator \cite{HL24}. An introduction of backstepping technique can be found in \cite{Krstic08}. Concerning the Neumann boundary control on the right, for a subclass of critical lengths, a time-varying feedback was given in \cite{CRX17} for which an exponential decay rate holds but cannot be arbitrary. It is interesting to know whether or not a similar phenomenon holds for \eqref{sys-KdV}. It is completely open to obtain the rapid stabilization of \eqref{sys-KdV} even for  time-varying feedbacks for critical lengths. 

\subsection{Organisation of the paper} \Cref{sect-Pre}, we establish some results used in the proof of the stablization. The rapid stabilization is studied in \Cref{sect-W-feedback} and the finite-time stabilization is investigated in \Cref{sect-D-feedback}.

\section{Preliminaries} \label{sect-Pre}

In this section, we first give the meaning of the weak solutions used in \Cref{thm-W} and \Cref{thm-D}. We then state and prove several well-posedness and stability results on the KdV equation. We finally establish the upper bound of $\| Q\| = \| Q \|_{\cL(L^2(0, L))}$ and $\| Q^{-1}\| = \| Q^{-1} \|_{\cL(L^2(0, L))}$ with $Q = Q(\lambda)$ being defined in \eqref{def-Q-KdV}, where the dependence on $\lambda$ is explicit. 

\medskip 
We begin with the following definition. 

\begin{definition} \label{def-KdV1}
Let $L>0$,  $T > 0$, $M \in \cL(L^2(0, L))$, $y_0  \in L^2(0, L)$, $f \in L^1((0, T); L^2(0, L))$, and $h \in L^2(0, T)$. A function $y \in X_T$ is a weak solution of the system 
\be \label{thm-KdV-sys-def}
\left\{\begin{array}{cl}
y_t  + y_x + y_{xxx} + M y = f &\mbox{ in } (0, T) \times (0, L), \\[6pt]
y(\cdot, 0) = y (\cdot, L) = 0, y_x(\cdot, L) - y_x(\cdot, 0) = h & \mbox{ in } (0, T), \\[6pt]
y(0, \cdot) = y_0 & \mbox{ in } (0, L), 
\end{array} 
\right. 
\ee
if 
\begin{multline}
\int_0^T \int_0^L \Big(f(t, x) - M y(t, \cdot) \Big) \varphi (t, x) \, dx \, dt + \int_0^L y_0 (x) \varphi(0, x) \, dx   + \int_0^T  h (t) \varphi_x (t, L) \, dt \\[6pt]
= - \int_0^T \int_0^L y (\varphi_t + \varphi_x + \varphi_{xxx}) \, dx \, dt  
\end{multline}
for all $\varphi \in C^3([0, T] \times [0, L])$ with $\varphi(T, \cdot) = 0$ and  $\varphi(\cdot, 0) = \varphi(\cdot, L) = \varphi_x(\cdot, L) - \varphi_x(\cdot, 0)=0$. 
\end{definition}

Concerning the nonlinear setting involving $(y, \ty)$, we use the following definition.  

\begin{definition} \label{def-KdV2} Let $L>0$, $T>0$, $M, \tM \in \cL(L^2(0, L))$, and let $y_0, \ty_0 \in L^2(0, L)$.  
A pair of functions $(y, \ty) \in X_T \times X_T$ is a weak solution of 
\be
\left\{\begin{array}{cl}
y_t  + y_x + y_{xxx} + y y_x = 0 &\mbox{ in } (0, T) \times (0, L), \\[6pt]
\ty_t + \ty_x + \ty_{xxx} + \tM \ty +  M y +   \ty y_x = 0 & \mbox{ in } (0, T) \times (0, L), \\[6pt]
y(\cdot, 0) = y (\cdot, L) = 0, y_x(\cdot, L) - y_x(\cdot, 0) = - \ty_x (\cdot, L) & \mbox{ in } (0, T), \\[6pt]
\ty(\cdot, 0) = \ty (\cdot, L) = 0, \ty_x(\cdot, L) - \ty_x(\cdot, 0) = 0 & \mbox{ in } (0, T), \\[6pt]
y(0) = y_0, \quad  \ty(0) = \ty_0 & \mbox{ in } (0, L).  
\end{array} 
\right. 
\ee
if, under the form of \eqref{thm-KdV-sys-def},  $y$ is the solution of the system with the internal source term $f = - y y_x$ and the boundary source term $h = - \ty_x (\cdot, L)$, and $\ty$ is the solution of the corresponding system with the internal source term $\tf = - (\tM \ty + M y +   \ty y_x) $ and the boundary source term $0$. 
\end{definition}

\begin{remark} \rm These definitions are compatible with the ones given in the semi-group language, see e.g., \cite[Section 3]{Ng-Riccati}: the weak solutions given here are also the weak solutions given in the semigroup terminology in \cite[Section 3]{Ng-Riccati}. 
\end{remark}

We next discuss the well-posedness and the stability of the weak solutions. The following result is on the linear setting given in \eqref{def-KdV1}. 

\begin{lemma} \label{lem-KdV1} Let $L > 0$, $0 < T < T_0$, $M \in \cL (L^2(0, L))$, and let $y_0 \in L^2(0, L)$, $f \in L^1((0, T); L^2(0, L))$, and $h \in L^2(0, T)$. There exists a unique weak solution $z \in X_T$ of the system 
\be \label{lem-KdV1-sys}
\left\{\begin{array}{c}
y_t + y_x + y_{xxx} + M y = f \quad \mbox{ in } (0, T) \times (0, L), \\[6pt]
y(\cdot, 0) = y (\cdot, L) = 0, y_x(\cdot, L) - y_x(\cdot, 0) = h \mbox{ in } (0, T), \\[6pt]
y(0, \cdot) = y_0 (\cdot)  \mbox{ in } (0, L).   
\end{array} 
\right. 
\ee
Moreover, 
\be \label{lem-KdV1-est}
\| y\|_{X_T} + \|y_x(L, \cdot) \|_{L^2(0, T)} \le C \Big(\| y_0\|_{L^2(0, T)} +  \| f\|_{L^1((0, T)); L^2(0, T))} + \| h \|_{L^2(0, T)} \Big), 
\ee
for some positive constant $C$ independent of $f$, $h$, $y_0$, and $T$. 
\end{lemma}

\begin{proof} We begin with the case $M 
\equiv 0$ as follows.  We first note that in the case $f \equiv 0$ and $y_0 \equiv 0$, we have, for $\xi \in \mC$,  
$$
i \xi \hy + \hy_x + \hy_{xxx} = 0 \mbox{ in } \mR \times (0, L),  
$$
where 
$$
\hy(\xi, x) = \frac{1}{\sqrt{2 \pi}} \int_0^\infty y(t, x) e^{- it \xi} \, dt. 
$$
For $\xi \in \mC$,  let $\lambda_j = \lambda_j(\xi)$ with $j=1, 2, 3$ be the three solutions of the equation $\lambda^3 + \lambda + i \xi = 0$. 
Taking into account the equation of $\hy$, we search for the solution of the form
\begin{equation}\label{lem-kdv2-p0}
\hy(\xi, \cdot) = \sum_{j=1}^3 a_j  e^{\lambda_j x},
\end{equation}
where  $a_j = a_j(\xi)$ for $j = 1, 2, 3$. Using the boundary condition, we then have 
\begin{equation*}
\begin{array}{c}
\sum_{j=1}^3 a_j  =0, \\[6pt]
\sum_{j=1}^3 a_j e^{\lambda_j L}  =0, \\[6pt]
\sum_{j=1}^3 a_j \lambda_j (e^{\lambda_j L} -1)  = \hh, 
\end{array}
\end{equation*}
where $$
\hh(\xi) = \frac{1}{\sqrt{2 \pi}} \int_0^\infty h(t) e^{- it \xi} \, dt. 
$$
This implies, with the convention $\lambda_{j+3} = \lambda_j$,  for $\xi \in \mR$, 
\begin{equation}\label{lem-kdv1-p1}
a_j =  \frac{e^{\lambda_{j+2}L} - e^{\lambda_{j+1} L}}{\det Q} \hh_3 \mbox{ for } j=1, 2, 3, 
\end{equation}
where 
$$
Q = Q (\xi) : = 
 \begin{pmatrix}
  1&1&1 \\
  e^{\lambda_1L}&e^{\lambda_2L}&e^{\lambda_3L}\\
  \lambda_1 (e^{\lambda_1  L} -1) & \lambda_2 (e^{\lambda_2  L} -1) & \lambda_3 (e^{\lambda_3  L} -1)
 \end{pmatrix}. 
$$
As in the proof of \cite[Lemma 4.4]{CKN-24} or \cite[Proposition 3.1]{Ng-KdV-D}, one can show that there exists a solution $y \in X_T$ satisfying 
\eqref{lem-KdV1-sys} and  \eqref{lem-KdV1-est} with $M \equiv 0$, $f\equiv0$, and $y_0 \equiv 0$. The existence of a  solution $y \in X_T$ satisfying 
\eqref{lem-KdV1-sys} and  \eqref{lem-KdV1-est} with $M \equiv 0$ for a general $f$, $y_0$, and $h$ follows from \cite[(4.17)]{Rosier97}. 

The proof of the uniqueness in the case $M \equiv 0$ can be proceeded as in \cite{Ng-KdV-D} (see also \cite[Chapter 8]{Coron07}). 
Let $y \in X_T$ be a solution with the zero data, i.e., $f \equiv 0$, $y_0 \equiv 0$, and $h \equiv 0$. Fix $\psi \in C^\infty_c \big((0, T) \times (0, L) \big)$ (arbitrarily). Let $\ty \in X_T$ be a solution of the backward system  
\begin{equation*}\left\{
\begin{array}{cl}
\ty_{t} + \ty_{x}  + \ty_{xxx}  = \psi &  \mbox{ in } (0, T) \times  (0, L),
\\[6pt]
\ty(\cdot, 0) = 0,  \;  \ty(\cdot, L) = 0, \;  \ty_x(\cdot, L) - \ty_x(\cdot, 0)   = 0 & \mbox{ in } (0,
T), \\[6pt]
\ty(T, \cdot)  = 0 & \mbox{ in } (0, L).
\end{array}\right.
\end{equation*}
Using the construction given previously, one can assume that $\ty$ is smooth. Using the definition of the weak solutions, we derive that 
$$
\int_{0}^T \int_0^L \psi (t, x) y(t, x) \, dt \, dx =0. 
$$
Since $\psi \in C^\infty_c \big((0, T) \times (0, L) \big)$ is arbitrary, we deduce that 
$$
y = 0 \mbox{ in } (0, T) \times (0, L). 
$$
The uniqueness is proved in the case $M \equiv 0$.

\medskip
The proof in the general case where $M$ is not required to be 0 follows from the case $M \equiv 0$ by using appropriate weighted norms involving time in $X_T$, see,  e.g., \cite[Section 4]{Ng-Riccati}. The details are omitted. 
\end{proof}

The following result on a specific linear setting is useful. 

\begin{lemma} \label{lem-KdV2} Let $L > 0$, $0 < T < T_0$, $\lambda \ge \lambda_0  > 0$,  and let $\ty_0 \in L^2(0, L)$ and $f \in L^1((0, T); L^2(0, L))$. Let $\ty \in X_T$ be the unique weak solution of the system 
\be
\left\{\begin{array}{cl}
\ty_t + \ty_x + \ty_{xxx} + 2 \lambda \ty  = \tf & \mbox{ in } (0, T) \times (0, L), \\[6pt]
\ty(\cdot, 0) = \ty (\cdot, L) = 0, \ty_x(\cdot, L) - \ty_x(\cdot, 0) = 0 & \mbox{ in } (0, T), \\[6pt]
\ty(0, \cdot) = \ty_0 & \mbox{ in } (0, L).  
\end{array} 
\right. 
\ee
Then 
\be
\| \ty(t, \cdot) \|_{L^2(0, T)} \le  e^{- 2 \lambda t}    \|\ty (0, \cdot) \|_{L^2(0, T)} + C \| \tf \|_{L^1((0, T)); L^2(0, L))} 
\ee 
and
\be
\| \ty \|_{X_T} + \| \ty_x(L, \cdot) \|_{L^2(0, T)} \le  C \Big( \|\ty (0, \cdot) \|_{L^2(0, T)} +  \lambda \| \tf \|_{L^1((0, T)); L^2(0, L))} \Big)
\ee
for some positive constant $C$ independent of $\tf$, $\lambda$, and $T$. 
\end{lemma}

\begin{remark} \rm The difference between \Cref{lem-KdV2} and \Cref{lem-KdV1} is the explicit dependence on the parameter $\lambda$ in \Cref{lem-KdV2}. This is useful to establish the finite-time stabilization result. 
\end{remark}

\begin{proof} Set 
$$
\ty_{2\lambda} (t, x) = \ty (t, x) e^{2 \lambda t} \mbox{ in } (0, T) \times (0, L). 
$$
We have 
\be \label{lem-KdV-2-p1}
\left\{\begin{array}{cl}
\ty_{2\lambda, t} + \ty_{2\lambda, x} + \ty_{2\lambda, xxx}  = e^{2 \lambda t} \tf  & \mbox{ in } (0, T) \times (0, L), \\[6pt]
\ty_{2\lambda}(\cdot, 0) = \ty_{2\lambda} (\cdot, L) = 0, \ty_{2\lambda, x} (\cdot, L) - \ty_{2\lambda, x}(\cdot, 0) = 0 & \mbox{ in } (0, T), \\[6pt]
\ty_{2\lambda} (0, \cdot) = \ty_0 & \mbox{ in } (0, L).  
\end{array} 
\right. 
\ee

By the linearity of the system, it suffices to consider two cases $\tf \equiv 0$ and $\ty_0 \equiv 0$ separately.

We first consider the case $\tf \equiv 0$.  Applying \Cref{lem-KdV1} to $\ty_{2 \lambda}$, we obtain 
\be \label{lem-KdV2-p1-1}
\| \ty_{2\lambda} \|_{X_T} + \| \ty_{2\lambda, x} (L, \cdot) \|_{L^2(0, T)} \le  C    \|\ty_{2\lambda}(0, \cdot) \|_{L^2(0, T)}, 
\ee 
and, since $A$ is skew-adjoint,  
\be \label{lem-KdV2-p1-2}
\| \ty_{2\lambda}(t, \cdot) \|_{L^2(0, T)} =   \|\ty_{2\lambda} (0, \cdot) \|_{L^2(0, T)}.
\ee
We derive from \eqref{lem-KdV2-p1-1} that 
\be
\| \ty(t, \cdot) \|_{L^2(0, T)} =  e^{- 2 \lambda t}    \|\ty (0, \cdot) \|_{L^2(0, T)}. 
\ee
Considering the system of $\ty$ and viewing $\- 2 \lambda \ty$ as a source term, after applying \Cref{lem-KdV1} to $\ty$, we obtain  
\be
\| \ty \|_{X_T} + \| \ty_{x}(L, \cdot) \|_{L^2(0, T)} \le  C   \|\ty(0, \cdot) \|_{L^2(0, T)}. 
\ee

We next consider the case $\ty_0 \equiv 0$. Applying \Cref{lem-KdV1}, we derive from \eqref{lem-KdV-2-p1} that, for $0< t \le T$, 
$$
\| \ty_{2\lambda} \|_{X_t} + \| \ty_{2 \lambda, x}(L, \cdot) \|_{L^2(0, t)} \le  C \| f e^{2 \lambda s}\|_{L^1((0, t)); L^2(0, L))}.  
$$ 
This implies 
\be \label{lem-KdV2-p2-1}
\| \ty(t, \cdot) \|_{L^2(0, L)}  \le C \| f \|_{L^1((0, T)); L^2(0, L))}.  
\ee 
Considering the system of $\ty$ and viewing $\tf - 2 \lambda \ty$ as a source, after applying \Cref{lem-KdV1} to $\ty$, we obtain 
\be \label{lem-KdV2-p2-2}
\| \ty \|_{X_T} + \| \ty_x(L, \cdot) \|_{L^2(0, T)} \le  C \lambda \| f \|_{L^1((0, T)); L^2(0, L))}.  
\ee 
The conclusions in the case $y_0 \equiv 0$ follow from \eqref{lem-KdV2-p2-1} and \eqref{lem-KdV2-p2-2}. 

\medskip
The proof is complete. 
\end{proof}

We next derive an upper bound and a lower bound for $\| Q\| $ with $Q = Q(\lambda)$ for which the dependence on $\lambda$ is explicit.  We begin with a result on an upper bound for $\|Q\|$, which is a consequence of the admissibility of the control \eqref{admissibility1}. 

\begin{lemma} \label{lem-Qlambda-U} Let $L \not \in \cN$, $\lambda \ge  \lambda_0 > 0$ and let $Q = Q(\lambda)$ be defined by \eqref{def-Q-KdV}. There exists a positive constant $C$ independent of $\lambda$ such that 
\be \label{lem-Qlambda-U-est}
\langle Q z, z \rangle_{L^2(0, L)} \le C \| z\|_{L^2(0, L)}^2 \mbox{ for all } z \in L^2(0, L). 
\ee
\end{lemma}

\begin{proof} We have
\begin{multline*}
\langle Q z, z \rangle_{L^2(0, L)} = \int_{0}^\infty e^{-2 \lambda s} | B^* e^{-s A^*} z|^2 \, ds = \sum_{n \ge 0} \int_{n}^{n+1} e^{-2 \lambda s} | B^* e^{-s A^*} z|^2 \, ds \\[6pt] 
\mathop{\le}^{\eqref{A-skew-adjoint}} \sum_{n \ge 0} e^{-2 \lambda n} \int_{n}^{n+1} | B^* e^{s A} z|^2 \, ds  \mathop{\le}^{\eqref{admissibility1-***}} \sum_{n \ge 0} C e^{-2 \lambda n} \| e^{n A} z \|_{L^2(0, L)}^2. 
\end{multline*}
The conclusion follows since $\| e^{n A} z \|_{L^2(0, L)} = \| z \|_{L^2(0, L)}$ thanks to the fact that $A$ is skew-adjoint. 
\end{proof}

\begin{remark} \rm As a consequence of \eqref{lem-Qlambda-U-est}, we derive that 
$$
\| Q (\lambda) \| \le C \mbox{ for } \lambda \ge \lambda_0, 
$$
for some positive constant $C$ independent of $\lambda$. 
\end{remark}

We next derive a lower bound for the norm $\|Q(\lambda)\|$ when $L \not \in \cN$, which implies an upper bound for $\|Q(\lambda)^{-1}\|$. To this end, we first state an observability inequality,  which is a consequence of a result of Lissy on the cost of controls for small time of the KdV equation \cite[Theorem 3.4]{Lissy-14}. 

\begin{proposition}\cite[Theorem 3.4]{Lissy-14} \label{pro-cost} Let $L \not \in \cN$ and $0 < T < T_0$.  We have, for some positive constant $C$ independent of $T$, 
$$
\int_0^T | B^* e^{-s A^*} z|^2  \ge e^{- \frac{C}{T^{1/2}}} \|z\|_{L^2(0, L)}^2 \mbox{ for all } z \in L^2(0, L).  
$$
\end{proposition}

Using \Cref{pro-cost}, we can prove the following result.  

\begin{lemma} \label{lem-Qlambda} Let $L \not \in \cN$, $\lambda \ge  \lambda_0 > 0$ and let $Q = Q(\lambda)$ be defined by \eqref{def-Q-KdV}. There exists a positive constant $C$ independent of $\lambda$ such that 
\be \label{lem-Qlambda-est}
\langle Q z, z \rangle_{L^2(0, L)} \ge e^{-C \lambda^{1/3}} \| z\|_{L^2(0, L)}^2 \mbox{ for all } z \in L^2(0, L). 
\ee
\end{lemma}

\begin{proof}
We have 
\begin{multline*}
\langle Q z, z \rangle_{L^2(0, L)} = \int_{0}^\infty e^{-2 \lambda s} | B^* e^{-s A^*} z|^2 \, ds \ge \int_{\lambda^{-2/3}}^{2\lambda^{-2/3}} e^{-2 \lambda s} | B^* e^{-s A^*} z|^2 \, ds \\[6pt]
\ge \lambda^{-2/3}  e^{-4 \lambda^{1/3} } \int_{\lambda^{-2/3}}^{2\lambda^{-2/3}}  | B^* e^{-s A^*} z|^2  
\mathop{\ge}^{\Cref{pro-cost}}  \lambda^{-2/3}  e^{-4 \lambda^{1/3} } e^{- C \lambda^{1/3}} \| e^{- \lambda^{-2/3} A^*} z\|_{L^2(0, L)}. 
\end{multline*}
This implies, since $A$ is skew-adjoint, 
$$
\langle Q z, z \rangle_{L^2(0, L)} \ge  e^{- C \lambda^{1/3}} \| z\|_{L^2(0, L)}, 
$$
which is the conclusion. 
\end{proof}

\begin{remark} \rm As a consequence of \eqref{lem-Qlambda-est}, we derive that 
$$
\| Q (\lambda)^{-1} \| \le  e^{C \lambda^{1/3}} \mbox{ for } \lambda \ge \lambda_0, 
$$
for some positive constant $C$ independent of $\lambda$. 
\end{remark}

\section{Static feedback in the trajectory sense for the KdV system} \label{sect-W-feedback}

This section is devoted to the proof of \Cref{thm-W} and its finite-time stabilization version \Cref{pro-FT}.  The key point of the analysis of this section is the following result which in particularly implies \Cref{thm-W}. This result is also the key ingredient of the proof of  \Cref{pro-FT} below on the local finite-time stabilization. 

\begin{proposition} \label{pro-KdV}
Let $L \not \in \cN$,  $0 < T \le T_0$, and $\lambda \ge  \lambda_0 > 0$, and let $Q = Q(\lambda)$ be defined by \eqref{def-Q-KdV}.  There exists $\eps > 0$ depending only on $L$, $T_0$, and $\lambda_0$ such that for all $y_0 \in L^2(0, L)$ if 
\be \label{lem-KdV-small}
(\lambda + e^{\frac{4}{3}\lambda T}) \Big( \| Q^{-1} y_0 \|_{L^2(0, L)} + \| y_0 \|_{L^2(0, L)} \Big) \le \eps,  
\ee
then there exists a unique weak solution $(y, \ty) \in X_T \times X_T$ of the following system 
\be \label{lem-KdV-sys}
\left\{\begin{array}{c}
y_t  + y_x + y_{xxx} + y y_x = 0 \quad \mbox{ in } (0, T), \\[6pt]
\ty_t + \ty_x + \ty_{xxx} + 2 \lambda \ty +   \ty y_x = 0 \mbox{ in } (0, T), \\[6pt]
y(\cdot, 0) = y (\cdot, L) = 0, y_x(\cdot, L) - y_x(\cdot, 0) = - \ty_x (\cdot, L) \mbox{ in } (0, T), \\[6pt]
\ty(\cdot, 0) = \ty (\cdot, L) = 0, \ty_x(\cdot, L) - \ty_x(\cdot, 0) = 0 \mbox{ in } (0, T), \\[6pt]
y(0) = y_0, \quad  \ty(0) = \ty_0: = Q^{-1} y_0. 
\end{array} 
\right. 
\ee
Moreover, we have 
\be \label{lem-KdV-cl0}
y (t, \cdot)= Q \ty (t, \cdot) \mbox{ for } t \in [0, T],  
\ee
\be \label{lem-KdV-cl1}
\| \ty_x(\cdot, L) \|_{L^2(0, T)}  + \| y\|_{X_T} + \| \ty \|_{X_T} \le C \Big(\| y_0 \|_{L^2(0, T)} + \| \ty_0 \|_{L^2(0, T)} \Big), 
\ee
and 
\be \label{lem-KdV-cl2}
\| \ty (t, \cdot) \|_{L^2(0, L)}  \le 2 e^{- 2 \lambda t} \Big(\| y_0 \|_{L^2(0, T)} + \| \ty_0 \|_{L^2(0, T)} \Big) \mbox{ for } t \in [0, T]. 
\ee
\end{proposition}

\begin{proof}  
Concerning the well-posedness and the stability of $(y, \ty)$, we only give the proof of \eqref{lem-KdV-cl1} and \eqref{lem-KdV-cl2}. The well-posedness can be proceeded by a standard fixed point argument, which involves tha analysis used in the proof of \eqref{lem-KdV-cl1} and \eqref{lem-KdV-cl2}, and omitted.  Applying \Cref{lem-KdV1} to $y$ and \Cref{lem-KdV2} to $\ty$, we have 
\be \label{lem-KdV-p1}
\| y\|_{X_T} \le C \Big(\| y_0 \|_{L^2(0, L)} + \| \ty_x(\cdot, L) \|_{L^2(0, T)} + \| y y_x \|_{L^1((0, T); L^2(0, L))} \Big)
\ee
and 
\be \label{lem-KdV-p2}
\| \ty_x(\cdot, L) \|_{L^2(0, T)} + \| \ty\|_{X_T} \le C \Big(\| \ty_0 \|_{L^2(0, T)} + \lambda \| \ty y_x \|_{L^1((0, T); L^2(0, L))} \Big). 
\ee
Here and in what follows, $C$ denotes a positive constant depending only on $T_0$, $L$, and $\lambda_0$. 

Combining \eqref{lem-KdV-p1} and \eqref{lem-KdV-p2} yields 
\begin{multline} \label{lem-KdV-p3}
\| y\|_{X_T} + \| \ty \|_{X_T} + \| \ty_x(\cdot, L) \|_{L^2(0, T)} \\[6pt]\le C \Big(\| y_0 \|_{L^2(0, L)} + \| \ty_0 \|_{L^2(0, L)} + \| y y_x \|_{L^1((0, T); L^2(0, L))} + \lambda \| \ty y_x \|_{L^1((0, T); L^2(0, L))} \Big).
\end{multline}
Since 
$$
\| y y_x \|_{L^1((0, T); L^2(0, L))} \le C \| y\|_{X_T}^2 \quad \mbox{ and } \quad \| \ty y_x \|_{L^1((0, T); L^2(0, L))} \le C \| \ty \|_{X_T} \| y\|_{X_T}, 
$$
it follows from \eqref{lem-KdV-p3} that 
\be \label{lem-KdV-p4} 
\| y\|_{X_T} + \| \ty \|_{X_T} + \| \ty_x(\cdot, L) \|_{L^2(0, T)}   \le C \Big(\| y_0 \|_{L^2(0, T)} + \| \ty_0 \|_{L^2(0, T)} \Big)
\ee
provided that 
\be \label{lem-KdV-cond1}
\lambda (\| y_0 \|_{L^2(0, L)} + \| \ty_0 \|_{L^2(0, T)}) \le \alpha, 
\ee
for some small positive constant $\alpha$ depending only on $T_0$, $L$, and $\lambda_0$. Condition \eqref{lem-KdV-cond1} is assumed from later on.  Applying \Cref{lem-KdV2} to $\ty$, we obtain 
\be \label{lem-KdV-p41}
\| \ty(t, \cdot) \|_{L^2(0, L)}  \le e^{- 2 \lambda t}    \|\ty (0, \cdot) \|_{L^2(0, T)} + C  \| \ty  y_x \|_{L^1((0, T)); L^2(0, L))},  
\ee
which yields, by \eqref{lem-KdV-p4},  
\be \label{lem-KdV-part1}
\| \ty(t, \cdot) \|_{L^2(0, L)}  \le e^{- 2 \lambda t}    \|\ty (0, \cdot) \|_{L^2(0, T)}  + C \Big(\| y_0 \|_{L^2(0, T)} +  \| \ty_0 \|_{L^2(0, T)} \Big)^2. 
\ee

We now improve \eqref{lem-KdV-part1}.
Set 
$$
\ty_{2 \lambda} (t, x) = e^{2 \lambda t} \ty (t, x) \mbox{ in } (0, T) \times (0, L). 
$$
Then 
$$
\ty_{2 \lambda, t} + \ty_{2 \lambda, x} + \ty_{2 \lambda, xxx} =  -  e^{2 \lambda t}  \ty y_x  \mbox{ in } (0, T) \times (0, L). 
$$
We get, for $0 \le t \le T$,  
\be \label{lem-KdV-p42}
\| \ty_{2 \lambda} (t, \cdot) \|_{L^2(0, L)} \le  \| \ty_{2 \lambda} (0, \cdot) \|_{L^2(0, L)} + C \| \ty y_x e^{2 \lambda s}\|_{L^1((0, t)); L^2(0, L))}.  
\ee

We next estimate the last term. Using the interpolation inequality given in \Cref{lem-interpolation} below, 
we obtain 
\begin{multline} \label{lem-KdV-p5}
\int_0^t \Big(\int_0^L |\ty (s, x) y_x (s, x) e^{2 \lambda s}|^2 \, dx \Big)^{1/2} \, ds \\[6pt]
\le \int_0^t e^{2 \lambda s}  \| \ty (s, \cdot) \|_{L^2(0, L)}^{1/2} \| \ty_x (s, \cdot) \|_{L^2(0, L)}^{1/2}
\| y_x (s, \cdot)\|_{L^2(0, L)}. 
\end{multline}
Using \eqref{lem-KdV-part1}, we derive from \eqref{lem-KdV-p5} that 
\begin{multline*}
\int_0^t \Big(\int_0^L |\ty (s, x) y_x (s, x) e^{2 \lambda s}|^2 \, dx \Big)^{1/2} \, ds \\[6pt]
\le C  \int_0^t e^{2 \lambda s}  \Big(e^{- 2 \lambda s}    \|\ty (0, \cdot) \|_{L^2(0, T)}  + \| y_0 \|_{L^2(0, T)}^2 +  \| \ty_0 \|_{L^2(0, T)}^2 \Big)^{1/2} \| \ty_x (s, \cdot) \|_{L^2(0, L)}^{1/2} \| y_x (s, \cdot)\|_{L^2(0, L)}.  
\end{multline*}
This yields, by \eqref{lem-KdV-p4}, 
\begin{multline} \label{lem-KdV-p6}
\int_0^t \Big(\int_0^L |\ty (s, x) y_x (s, x) e^{2 \lambda s}|^2 \, dx \Big)^{1/2} \, ds
\le C e^{\lambda t}\Big(\| y_0 \|_{L^2(0, T)} +  \| \ty_0 \|_{L^2(0, T)} \Big)^{2}  \\[6pt] + C e^{2 \lambda t}\Big(\| y_0 \|_{L^2(0, T)} +  \| \ty_0 \|_{L^2(0, T)} \Big)^{5/2}. 
\end{multline}
Since $y_{2 \lambda} (t, x) = e^{2 \lambda t} \ty(t, x)$, we derive from \eqref{lem-KdV-p42} and \eqref{lem-KdV-p6} that 
\begin{multline} \label{lem-KdV-p6-1}
\| \ty (t, \cdot) \|_{L^2(0, L)}  \le e^{- 2 \lambda t} \| \ty_0 \|_{L^2(0, L)} + C e^{- \lambda t} \Big(\| y_0 \|_{L^2(0, T)} + \| \ty_0 \|_{L^2(0, T)} \Big)^{2} \\[6pt]  + C \Big(\| y_0 \|_{L^2(0, T)} + \| \ty_0 \|_{L^2(0, T)} \Big)^{5/2}.  
\end{multline}
We deduce from \eqref{lem-KdV-p6-1} that  
\be \label{lem-KdV-part2}
\| \ty (t, \cdot) \|_{L^2(0, L)}  \le 2 e^{- 2 \lambda t} \Big( \| \ty_0 \|_{L^2(0, L)} + \| y_0 \|_{L^2(0, L)} \Big)
\ee
as long as 
\be \label{lem-KdV-cond2}
e^{\frac{4}{3} \lambda T} \Big( \| \ty_0 \|_{L^2(0, L)} + \| y_0 \|_{L^2(0, L)} \Big) \le \alpha 
\ee
for some positive constant $\alpha$ depending only on $T_0$, $L$, and $\lambda_0$. 

\medskip 
Assertion \eqref{lem-KdV-cl1} and \eqref{lem-KdV-cl2} now follows from \eqref{lem-KdV-p4} and \eqref{lem-KdV-part2} after noting \eqref{lem-KdV-cond1} and \eqref{lem-KdV-cond2}. 

\medskip 
We next establish \eqref{lem-KdV-cl0} in the spirit of \cite{Ng-Riccati}, which gives the meaning of the feedback in the {\it trajectory} sense. 
Set, for $t  \in [0, T]$ and $x \in (0, L)$, with  $f(t, x) = - y (t, x) y_x (t, x)$ and $\tf (t, x) = - \ty (t, x) y_x (t, x)$, 
\be
y_\lambda (t, x) = e^{\lambda t} y (t, x),  \quad \ty_\lambda (t, x) = e^{\lambda t} \ty (t, x), 
\ee
\be
f_\lambda (t, x) = e^{\lambda t} f (t, x), \quad \mbox{ and } \quad \tf_\lambda (t, x) = e^{\lambda t} \tf (t, x), 
\ee
and denote
$$
A_\lambda = A + \lambda I. 
$$
We have, since $A^* = - A$, 
\be \label{lem-KdV-syslambda}
\left\{\begin{array}{c}
y_\lambda' = A_\lambda y_\lambda - B B^* \ty_\lambda  + f_
\lambda \quad \mbox{ in } (0, T), \\[6pt]
\ty_\lambda'  = - A_\lambda^* \ty_\lambda  + \tf_\lambda   \quad \mbox{ in } (0, T), \\[6pt]
y_\lambda(0) = y(0), \quad  \ty_\lambda(0) = \ty(0). 
\end{array} 
\right. 
\ee

Set, for $t  \ge 0$, 
$$
 Z_\lambda(t) = y_\lambda (t) - Q \ty_\lambda (t).
$$

We formally have  
\begin{multline*} 
\frac{d}{dt} Z_\lambda = A_\lambda y_\lambda - B  B^* \ty_\lambda + f_\lambda + Q  A_\lambda^* \ty_\lambda  - Q \tf_{\lambda} \\[6pt]
= A_\lambda (y_\lambda  - Q \ty_\lambda) + A_\lambda Q  \ty_\lambda  - B  B^* \ty_\lambda+ Q  A_\lambda^* \ty_\lambda  - Q \tf_\lambda, 
\end{multline*}
which yields, by \eqref{identity-Op},   
\be\label{lem-KdV-cl0-p1}
\frac{d}{dt} Z_\lambda = A_\lambda Z_\lambda  + f_\lambda - Q \tf_\lambda. 
\ee

\medskip 
We now give the proof of \eqref{lem-KdV-cl0-p1} using the results in \cite{Ng-Riccati}. 
Let  $\tau > 0$,  $\varphi_\tau \in L^2(0, L)$ and let $\varphi \in C([0, \tau]; L^2(0, L))$ be the unique weak solution of the system 
\be \label{lem-KdV-varphi} 
\left\{\begin{array}{c}
\varphi ' = - A_\lambda^* \varphi  \mbox{ in } (0, \tau), \\[6pt]
\varphi(\tau) = \varphi_\tau  
\end{array} \right. 
\ee
(see, e.g., \cite[Section 3]{Ng-Riccati} for the definition of the weak solutions for which $\varphi$ satisfies).

Applying \cite[Lemma 3.1]{Ng-Riccati} for $A_\lambda$ with $t = \tau$, we derive from \eqref{lem-KdV-sys} and \eqref{lem-KdV-varphi} that  
\be\label{lem-KdV-cl0-p2}
 \langle y_\lambda(\tau), \varphi(\tau) \rangle_{L^2(0, L)} - \langle y_\lambda (0), \varphi(0) \rangle_{L^2(0, L)} =  - \int_0^\tau \langle  B^* \ty_\lambda (s), B^* \varphi(s) \rangle \, ds + \int_0^\tau \langle f_\lambda(s, \cdot),  \varphi(s, \cdot) \rangle \, ds.
\ee
 Using \eqref{identity-Op} and applying \cite[Lemma 4.1]{Ng-Riccati} for $A_\lambda$, $\ty_\lambda(\tau - \cdot)$ and $\varphi (\tau - \cdot)$, we obtain   
\begin{multline}\label{lem-KdV-cl0-p3}
\langle Q \ty_\lambda(0), \varphi(0) \rangle - \langle Q \ty_\lambda(\tau), \varphi (\tau) \rangle  \\[6pt]
 =  \int_0^\tau  \langle B^*\ty_\lambda(\tau - s), B^*\varphi (\tau - s) \rangle  \, d s  -  \int_0^\tau \langle Q \tf_\lambda(\tau -s, \cdot),  \varphi(\tau -s, \cdot) \rangle \, ds.   
\end{multline}
Summing \eqref{lem-KdV-cl0-p2} and \eqref{lem-KdV-cl0-p3},   we deduce from \eqref{lem-KdV-syslambda} and \eqref{lem-KdV-varphi} that  
$$
\langle Z_\lambda(\tau), \varphi(\tau) \rangle - \langle Z_\lambda(0), \varphi(0) \rangle = \int_0^\tau \langle f_\lambda(s, \cdot) - Q \tf_\lambda(s, \cdot),  \varphi(s, \cdot)  \rangle \, ds. 
$$
This yields 
$$
\langle Z_\lambda(\tau), \varphi(\tau) \rangle - \langle Z_\lambda(0), e^{\tau A^*}\varphi(\tau) \rangle  
=  \int_0^\tau \langle f_\lambda(s, \cdot) - Q \tf_\lambda(s, \cdot),  e^{(\tau -s) A^*}\varphi(\tau, \cdot) \rangle \, ds. 
$$
Since $\varphi(\tau) \in \mH$ is arbitrary, we obtain
$$
Z_\lambda(\tau) = e^{\tau A} Z_\lambda(0) + \int_0^\tau  e^{(\tau -s) A} \big( f_\lambda(s, \cdot) - Q \tf_\lambda(s, \cdot) 
\big)
, ds 
$$
which implies \eqref{lem-KdV-cl0-p1} (see also \cite[Section 3]{Ng-Riccati}). 

Note that 
$$
f_\lambda - Q \tf_\lambda = e^{\lambda t} (y- Q \ty) y_x = Z_
\lambda y_x. 
$$
Assertion \eqref{lem-KdV-cl0} follows from \eqref{lem-KdV-cl0-p1} for $\eps$ sufficiently small. 
\medskip 
The proof is complete. 
\end{proof}

In the proof of \Cref{pro-KdV}, we used the following simple interpolation inequality. 
\begin{lemma}\label{lem-interpolation} Let $L>0$. We have
$$
\| \varphi \|_{L^\infty(0, L)} \le  \| \varphi \|_{L^2(0, L)}^{1/2}\| \varphi' \|_{L^2(0, L)}^{1/2} \mbox{ for } \varphi \in H^1(0, L) \mbox{ with } \varphi(0) = 0.  
$$
\end{lemma}

\begin{proof} The result is just a consequence of the fact, for $x \in [0, L]$,  
$$
\varphi^2(x) = \varphi^2(x) - \varphi^2(0) = 2 \int_0^x \varphi'(s) \varphi(s) \, ds  
$$
and the H\"older inequality. 
\end{proof}

We next state and prove the finite-time stabilization result in the trajectory sense.

\begin{proposition} \label{pro-FT}  Let $L \not \in \cN$ and $T > 0$.  Let $(t_n)$ be an increasing sequence that converges to $T$ with $t_0 = 0$ and let $(\lambda_n) \subset \mR_+$ be an increasing sequence.  Set $s_0 = 0$ and $s_n = \sum_{k=0}^{n-1} \lambda_k(t_{k+1} - t_k)$ for $n \ge 1$.  There exists a constant $\gamma > 1$ such that, if for large $n$, 
\be \label{pro-FT-lambda1}
(t_{n+1} - t_n) \lambda_n \ge \gamma \lambda_n^{1/3},  \quad \mbox{ and }  \quad 
\lambda_{n+1} (t_{n+2} - t_{n+1}) \le (1 + 1/\gamma) \lambda_n (t_{n+1} - t_n), 
\ee
and 
\be \label{pro-FT-lambda2}
\lim_{n \to + \infty} \frac{s_n}{n + \lambda_{n+1}^{1/3}} = + \infty, 
\ee
then there exists $\eps_0 > 0$ such that for all $y_0 \in L^2(0, L)$ with $\| y_0 \|_{L^2(0, L)} \le \eps_0$, there exists a unique pair $(y, \ty)$ be such that 
$y \in X_T$ and $\ty \in C([t_n, t_{n+1}); L^2(0, L)) \cap L^2( (t_n, t_{n+1}); H^1(0, L))$, 
for $t_{n} \le t < t_{n+1}$ and $n \ge 0$, and, for $n 
\ge 0$, it holds
\be \label{cor-KdV-sys}
\left\{\begin{array}{cl}
y_t + y_x + y_{xxx} + y y_x = 0 &  \mbox{ in } (t_n, t_{n+1}), \\[6pt]
\ty_t + \ty_x + \ty_{xxx} + 2 \lambda \ty +   \ty y_x = 0 & \mbox{ in } (t_n, t_{n+1}), \\[6pt]
y(\cdot, 0) = y (\cdot, L) = 0, y_x(\cdot, L) - y_x(\cdot, 0) = - \ty_x (\cdot, L) & \mbox{ in } (t_n, t_{n+1}), \\[6pt]
\ty(\cdot, 0) = \ty (\cdot, L) = 0, \ty_x(\cdot, L) - \ty_x(\cdot, 0) = 0 & \mbox{ in } (t_n, t_{n+1}), \\[6pt]
y(t_n, \cdot) = y(t_n, \cdot), \quad  \ty(t_n, \cdot) = Q_n^{-1} y(t_n, \cdot) & \mbox{ in } (0, L), 
\end{array} 
\right. 
\ee
where $Q_n = Q(\lambda_n)$ defined by \eqref{def-Q-KdV-PDE} with $\lambda = \lambda_n$. Moreover, we have, for $t_{n-1} \le t \le t_{n}$ and for $n \ge 1$,  
$$
\|y(t, \cdot) \|_{L^2(0, L)} \le e^{ - s_{n-1} + C (n + \lambda_{n-1}^{1/3}}) \| y_0\|_{L^2(0, L)}. 
$$
for some positive constant $C$ independent of $n$ and $y_0$. Consequently, it holds 
$$
y(t, \cdot) \to 0 \mbox{ in } L^2(0, L) \mbox{ as } t \to T_{-}. 
$$
\end{proposition}

\begin{remark} \rm There are sequences $(t_n)$ and $(\lambda_n)$ which satisfy the conditions given in the above proposition, for example,  $t_n = T - T/n^2$ and $\lambda_n = n^8$ for large $n$. 
\end{remark}

\begin{proof} Applying \Cref{pro-KdV} and  \Cref{lem-Qlambda-U,lem-Qlambda}, we have 
\be \label{pro-FT-p1}
\| y(t_n, \cdot) \|_{L^2(0, L)} \le  e^{-  2 \lambda_{n-1} (t_n - t_{n-1}) + C (1+ \lambda_{n-1}^{1/3})} \| y(t_{n-1}, \cdot) \|_{L^2(0, L)} \mbox{ for } n \ge 1.
\ee
It follows that, for $\gamma$ sufficiently large,  
\be
\| y(t_{n}, \cdot) \|_{L^2(0, L)} \le e^{- s_{n-1} + C n} \| y_0 \|_{L^2(0, L)} \mbox{ for } n \ge 1. 
\ee
The conclusion now follows from \Cref{pro-KdV}. 

It is worth noting that \eqref{pro-FT-lambda1} assure the existence of $(y, \ty)$ by applying \Cref{pro-KdV} in the time interval $(t_n, t_{n+1})$ for all $n \ge 1$ (after a translation of time) since the condition \eqref{lem-KdV-small} corresponding to the time interval $(t_n, t_{n+1})$, i.e., 
$$
(\lambda_n + e^{\frac{4}{3}\lambda_n (t_{n+1} - t_n)}) \Big( \| Q_n^{-1} y (t_n, \cdot) \|_{L^2(0, L)} + \| y (t_n, \cdot) \|_{L^2(0, L)} \Big) \le \eps,  
$$
is ensured if the following condition holds, for large $n$,  
$$
e^{\frac{4}{3} \lambda_n (t_{n+1} - t_n) + C \lambda_n^{1/3}} \| y(t_n, \cdot)\|_{L^2(0, L)} \le \eps. 
$$
if $\eps_0$ is sufficiently small. This holds by \eqref{pro-FT-lambda1} and \eqref{pro-FT-p1} if $\gamma$ is sufficiently large. 
\end{proof}

\section{Dynamic feedback for the KdV system} \label{sect-D-feedback}

This section is devoted to the proof of \Cref{thm-D} and its finite-time stabilization version \Cref{pro-FTD}.  The key point of the analysis of this section is the following result which not only implies \Cref{thm-D} but also is the key ingredient of the proof of  \Cref{pro-FTD}.

\begin{proposition}\label{pro-KdV-D} Let $L \not \in \cN$, $0< T \le T_0$,  $\lambda > \lambda_0 > 0$, $\lambda_1 > 0$, and $c_0 > 0$, and let $Q = Q(\lambda)$ be defined by \eqref{def-Q-KdV}. Assume that 
\be \label{lem-D-lambda}
\lambda_1 >  (2 + c_0) \lambda. 
\ee
There exists $\eps > 0$ depending only on $L$, $T_0$, $\lambda_0$, and $c_0$ such that for all $y_0, \ty_0 \in L^2(0, L)$ with  
\be\label{lem-KdVD-small}
\lambda^{2} e^{\frac{4}{3}\lambda T} \| Q^{-1} \|^{2} \Big( \| y_0 \|_{L^2(0, L)} + \| \ty_0 \|_{L^2(0, L)} \Big) \le \eps,  
\ee
then there exists a unique weak solution $(y, \ty) \in X_T \times X_T$ of the following system
\be \label{lem-D-sys}
\left\{\begin{array}{cl}
y_t  + y_x + y_{xxx} + y y_x =  0 &\mbox{ in } (0, T) \times (0, L), \\[6pt]
\ty_t + \ty_x + \ty_{xxx} + 2 \lambda \ty - 
\lambda_1 Q^{-1}(y - Q \ty) + \ty y_x = 0 & \mbox{ in } (0, T) \times (0, L), \\[6pt]
y(\cdot, 0) = y (\cdot, L) = 0, y_x(\cdot, L) - y_x(\cdot, 0) = - \ty_x (\cdot, L) & \mbox{ in } (0, T), \\[6pt]
\ty(\cdot, 0) = \ty (\cdot, L) = 0, \ty_x(\cdot, L) - \ty_x(\cdot, 0) = 0 & \mbox{ in } (0, T), \\[6pt]
y(0, \cdot) = y_0, \quad  \ty(0, \cdot) = \ty_0 & \mbox{ in } (0, L).  
\end{array} 
\right. 
\ee
Moreover, it holds 
\be \label{lem-D-cl1}
\| y\|_{X_T} + \| \ty\|_{X_T} \le C \lambda \| Q^{-1}\| \big( \| y(0, \cdot) \|_{L^2(0, L)} + \|\ty(0, \cdot)\|_{L^2(0, L)} \big) 
\ee
and
\be \label{lem-D-cl2}
\| y(t, \cdot)\|_{L^2(0, L)} + \| \ty(t, \cdot)\|_{L^2(0, L)} \le C \| Q^{-1} \|  e^{ - 2 \lambda t}  \big( \| y(0,  \cdot) \|_{L^2(0, L)} + \|\ty(0, \cdot)\|_{L^2(0, L)} \big) \mbox{ for } t \in [0, T], 
\ee
where $C$ is a positive constant independent of $\lambda$, $\lambda_1$,  $T$, and $(y_0, \ty_0)$. 
\end{proposition}

\begin{proof} Concerning the well-posedness and the stability of $(y, \ty)$, we only give the proof of \eqref{lem-D-cl1} and \eqref{lem-D-cl2}. The well-posedness can be proceeded by a standard fixed point argument and omitted.

Set, for $t \in [0, T]$ and $x \in (0, L)$, with  $f = - y y_x$ and $\tf = - \ty y_x$, 
\be
y_\lambda (t, x) = e^{\lambda t} y (t, x),  \quad \ty_\lambda (t, x) = e^{\lambda t} \ty (t, x), 
\ee
\be
f_\lambda (t, x) = e^{\lambda t} f (t, x),  \quad \mbox{ and } \quad \tf_\lambda (t, x) = e^{\lambda t} \tf (t, x), 
\ee
and denote
$$
A_\lambda = A + \lambda I. 
$$
We have
\be \label{lem-D-syslambda}
\left\{\begin{array}{c}
y_\lambda' = A_\lambda y_\lambda - B B^* \ty_\lambda  + f_
\lambda \quad \mbox{ in } (0, T), \\[6pt]
\ty_\lambda'  = - A_\lambda^* \ty_\lambda+ \lambda_1 Q^{-1} (y_\lambda - Q \ty_\lambda) + \tf_\lambda   \quad \mbox{ in } (0, T), \\[6pt]
y_\lambda(0) = y(0), \quad  \ty_\lambda(0) = \ty(0). 
\end{array} 
\right. 
\ee

Set, for $t  \ge 0$, 
$$
 Z_\lambda(t) = y_\lambda (t) - Q \ty_\lambda (t).
$$ 
As in the proof of \eqref{lem-KdV-cl0-p1} in the proof of \Cref{pro-KdV}, we have 
\be\label{lem-D-p1}
\frac{d}{dt} Z_\lambda = A_\lambda Z_\lambda - \lambda_1 Z_\lambda + f_\lambda - Q \tf_\lambda. 
\ee

We derive from \eqref{lem-D-p1} that, since $A$ is skew-adjoint, 
\be \label{lem-D-p00}
\| Z_\lambda (t, \cdot) \|_{L^2(0, L)} \le  e^{(-\lambda_1 + \lambda)t} \|Z_\lambda(0, \cdot) \|_{L^2(0, L)} + C e^{(-\lambda_1 + \lambda)t}  \| f  - Q \tf\|_{L^1((0, t); L^2(0, L))},    
\ee
which yields 
\be \label{lem-D-p11}
\| y(t, \cdot) - Q \ty(t, \cdot) \|_{L^2(0, L)} \le  e^{-\lambda_1 t} \| y(0, \cdot) - Q \ty(0, \cdot) \|_{L^2(0, L)} + C e^{-\lambda_1 t} \| f - Q \tf \|_{L^1((0, t); L^2(0, L))}.  
\ee
Here and in what follows in this proof, $C$ is a positive constant independent of $\lambda$, $t$, $T$, and $(y_0, \ty_0)$. 

Applying \Cref{lem-KdV1} to $y$ and \Cref{lem-KdV2} to $\ty$, we have 
\be \label{lem-D-c1}
\| y\|_{X_T} \le C \Big(\| y_0 \|_{L^2(0, L)} + \| \ty_x(\cdot, L) \|_{L^2(0, T)} + \| y y_x \|_{L^1((0, T); L^2(0, L))} \Big)
\ee
and 
\begin{multline} \label{lem-D-c2}
\| \ty_x(\cdot, L) \|_{L^2(0, T)} + \| \ty\|_{X_T} \\[6pt]
\le C \Big(\| \ty_0 \|_{L^2(0, T)} + \lambda \lambda_1 \| Q^{-1}\|  \| y - Q \ty \|_{L^1((0, T); L^2(0, L))} + 
\lambda \| \ty y_x \|_{L^1((0, T); L^2(0, L))} \Big). 
\end{multline}
Combining \eqref{lem-D-p11},  \eqref{lem-D-c1},  and \eqref{lem-D-c2} yield, since $\| Q^{-1}\| \ge C$ by \Cref{lem-Qlambda-U},  
\begin{multline} \label{lem-D-c3}
\| y\|_{X_T} + \| \ty \|_{X_T} + \| \ty_x(\cdot, L) \|_{L^2(0, T)} \\[6pt]\le C \lambda \| Q^{-1} \|\Big(\| y_0 \|_{L^2(0, L)} + \| \ty_0 \|_{L^2(0, L)} + \| y y_x \|_{L^1((0, T); L^2(0, L))} + \| \ty y_x \|_{L^1((0, T); L^2(0, L))} \Big).
\end{multline}
Since 
\be \label{lem-D-c4*} 
\| y y_x \|_{L^1((0, T); L^2(0, L))} \le C \| y\|_{X_T}^2 \quad \mbox{ and } \quad \| \ty y_x \|_{L^1((0, T); L^2(0, L))} \le C \| \ty \|_{X_T} \| y\|_{X_T}, 
\ee
it follows from \eqref{lem-D-c3} that 
\be \label{lem-D-c4} 
\| y\|_{X_T} + \| \ty \|_{X_T} + \| \ty_x(\cdot, L) \|_{L^2(0, T)}   \le C \lambda \| Q^{-1}\| \Big(\| y_0 \|_{L^2(0, T)} + \| \ty_0 \|_{L^2(0, T)} \Big)
\ee
provided that 
\be \label{lem-D-cond1}
\lambda \| Q^{-1}\| (\| y_0 \|_{L^2(0, L)} + \| \ty_0 \|_{L^2(0, T)}) \le \alpha, 
\ee
for some small positive constant $\alpha$ depending only on $T_0$, $L$, and $\lambda_0$; this condition is assumed from later on.  

Since 
$$
\ty'  = - A^* \ty - 2 \lambda \ty + \lambda_1 Q^{-1} (y - Q\ty) +   \tf\quad \mbox{ in } (0, T), 
$$
it follows that 
\be \label{lem-D-tylambda}
\ty_{2 \lambda}' = - A^* \ty_{2 \lambda} + \tg \mbox{ in } (0, T), 
\ee
where, for $t \in [0, T]$ and for $x \in (0, L)$, 
$$
\ty_{2 \lambda}(t, x) = e^{2 \lambda t} \ty (t, x) \quad \mbox{ and } \quad \tg(t, x) = \lambda_1 e^{2 \lambda t} Q^{-1} \big(y(t, x) - Q \ty (t, x)\big) + e^{2 \lambda t} \tf (t, x). 
$$
We thus obtain, since $A$ is skew-adjoint,  
\be\label{lem-D-p4}
\|\ty_{2 \lambda}(t, \cdot)\|_{L^2(0, L)} \le  \| \ty_{2 \lambda}(0, \cdot)\|_{L^2(0, L)} + C \| \tg \|_{L^1((0, t); L^2(0, L))}. 
\ee
Combining \eqref{lem-D-p11} and \eqref{lem-D-p4} yields 
\begin{multline*}
\| \ty_{2 \lambda} (t, \cdot) \|_{L^2(0, L)} \le C  \Big( \| \ty(0, \cdot) \|_{L^2(0, L)} + \| Q^{-1}\| \|y_0 - Q \ty_0 \|_{L^2(0, L)} \\[6pt]
+ \| Q^{-1}\| \big(\|f\|_{L^1((0, T); L^2(0, L))} + \|\tf\|_{L^1((0, T); L^2(0, L))} \big) + e^{2 \lambda t} \|\tf\|_{L^1((0, T); L^2(0, L))} \Big), 
\end{multline*}
which implies, since $\| Q^{-1}\| \ge C$ by \Cref{lem-Qlambda-U},  
\begin{multline}\label{lem-D-p22}
\| \ty(t, \cdot)\|_{L^2(0, L)} \le C e^{-2 \lambda t}\| Q^{-1}\| \big( \| y_0 \|_{L^2(0, L)} + \|\ty_0 \|_{L^2(0, L)} \big) \\[6pt]
+  C \| Q^{-1}\|  \big(  \|f\|_{L^1((0, T); L^2(0, L))} + \|\tf\|_{L^1((0, T); L^2(0, L))} \big). 
\end{multline}

Since $\| Q\| \le C$ by \Cref{lem-Qlambda-U}, we derive from \eqref{lem-D-p11}, \eqref{lem-D-c4},  and \eqref{lem-D-p22} that 
\begin{multline}\label{lem-D-p33}
\| y(t, \cdot)\|_{L^2(0, L)} + \| \ty(t, \cdot)\|_{L^2(0, L)} \le C e^{-2 \lambda t}\| Q^{-1}\| \big( \| y(0) \|_{L^2(0, L)} + \|\ty(0) \|_{L^2(0, L)} \big)  \\[6pt] 
+ C \lambda^2 \| Q^{-1} \|^{3}  \big( \| y(0) \|_{L^2(0, L)}^{2} + \|\ty(0) \|_{L^2(0, L)}^{2} \big).  
\end{multline}

We now improve \eqref{lem-D-p33}. From \eqref{lem-D-tylambda}, we obtain 
\begin{multline*}
\|\ty_{2 \lambda} (t, \cdot) \|_{L^2(0, L)} \\[6pt]
\le \|\ty_{2 \lambda} (0, \cdot) \|_{L^2(0, L)} + C \| \lambda_1 e^{2 \lambda s}Q^{-1} (y - Q \ty) \|_{L^1((0, t); L^2(0, L))} + C \|e^{2 \lambda s} \ty y_x \|_{L^1((0, t); L^2(0, L))}, 
\end{multline*} 
which yields, by \eqref{lem-D-p11} and \eqref{lem-D-c4}, and \Cref{lem-Qlambda-U},  
\begin{multline} \label{lem-D-t0}
\|\ty_{2 \lambda} (t, \cdot) \|_{L^2(0, L)}
\le \|\ty_{2 \lambda} (0, \cdot) \|_{L^2(0, L)} + C \| Q^{-1}\| \| y(0, \cdot) - Q \ty(0, \cdot) \|_{L^2(0, L)} \\[6pt] + C  \lambda^2 \|Q^{-1}\|^3 \Big(\| y_0 \|_{L^2(0, T)} + \| \ty_0 \|_{L^2(0, T)} \Big)^2 + C \|e^{2 \lambda s} \ty y_x \|_{L^1((0, t); L^2(0, L))}. 
\end{multline}
We have, by \Cref{lem-interpolation}, 
\begin{multline} \label{lem-D-t1}
\int_0^t \Big(\int_0^L |\ty (s, x) y_x (s, x) e^{2 \lambda s}|^2 \, dx \Big)^{1/2} \, ds \\[6pt]
\le \int_0^t e^{2 \lambda s}  \| \ty (s, \cdot) \|_{L^2(0, L)}^{1/2} \| \ty_x (s, \cdot) \|_{L^2(0, L)}^{1/2}
\| y_x (s, \cdot)\|_{L^2(0, L)}. 
\end{multline}
Using \eqref{lem-D-p33}, we derive from \eqref{lem-D-t1} that 
\begin{multline*}
\int_0^t \Big(\int_0^L |\ty (s, x) y_x (s, x) e^{2 \lambda s}|^2 \, dx \Big)^{1/2} \, ds 
\le C  \int_0^t e^{2 \lambda s}  \Big(e^{-2 \lambda s}\| Q^{-1}\| \big( \| y(0) \|_{L^2(0, L)} + \|\ty(0) \|_{L^2(0, L)} \big) \\[6pt]+  \lambda^2 \| Q^{-1}\|^{3} \big( \| y(0) \|_{L^2(0, L)} + \|\ty(0) \|_{L^2(0, L)} \big)^2 \Big)^{1/2} 
\| \ty_x (s, \cdot) \|_{L^2(0, L)}^{1/2} \| y_x (s, \cdot)\|_{L^2(0, L)} \, ds.  
\end{multline*}
This yields, by  \eqref{lem-D-c4}, 
\begin{multline}\label{lem-D-t01}
\int_0^t \Big(\int_0^L |\ty (s, x) y_x (s, x) e^{2 \lambda s}|^2 \, dx \Big)^{1/2} \, ds 
\le   C \lambda^{\frac{1}{2}} e^{\lambda t} \| Q^{-1}\|^{2} \big( \| y(0) \|_{L^2(0, L)}^{2} + \|\ty(0) \|_{L^2(0, L)}^{2} \big)   \\[6pt] 
 + C \lambda^{\frac{3}{2}} e^{2 \lambda t} \| Q^{-1}\|_{\cL(L^2(0, L))}^{3} \Big(\| y_0 \|_{L^2(0, T)}^\frac{5}{2} + \| \ty_0 \|_{L^2(0, T)}^\frac{5}{2} \Big).  
\end{multline}
Since $y_{2\lambda} (t, x) = e^{2 \lambda t} \ty(t, x)$, we derive from \eqref{lem-D-t0} and \eqref{lem-D-t01} that 
\begin{multline}\label{lem-D-t02}
\| \ty (t, \cdot) \|_{L^2(0, L)}  \le e^{- 2 \lambda t} \| \ty_0 \|_{L^2(0, L)} + 
C e^{- 2 \lambda t} \| Q^{-1}\| \| y(0, \cdot) - Q \ty(0, \cdot) \|_{L^2(0, L)}  \\[6pt] 
+ e^{-2 \lambda t } \lambda^2 \|Q^{-1}\|^3 \Big(\| y_0 \|_{L^2(0, T)}^2 + \| \ty_0 \|_{L^2(0, T)}^2 \Big)
+ C \lambda^\frac{1}{2} e^{- \lambda t } \|Q^{-1}\|^2 (\| y_0 \|_{L^2(0, T)}^2 + \| \ty_0 \|_{L^2(0, T)}^2 ) \\[6pt]
+ C \lambda^{\frac{3}{2}} \| Q^{-1}\|^{3} \Big(\| y_0 \|_{L^2(0, T)}^\frac{5}{2} + \| \ty_0 \|_{L^2(0, T)}^\frac{5}{2} \Big) . 
\end{multline}

We thus reach from \eqref{lem-D-t02} that 
\be\label{lem-D-t03}
\| \ty (t, \cdot) \|_{L^2(0, L)}  \le C \| Q\|^{-1} e^{- 2 \lambda t} \Big( \| \ty_0 \|_{L^2(0, L)} + \| y_0 \|_{L^2(0, L)} \Big), 
\ee
as long as 
\be\label{lem-D-cond2}
\lambda^2 e^{\frac{4}{3}\lambda T} \| Q^{-1} \|^{2} \Big( \| \ty_0 \|_{L^2(0, L)} + \| y_0 \|_{L^2(0, L)} \Big) \le \alpha. 
\ee

The conclusion now follows from \eqref{lem-D-p11} and \eqref{lem-D-t03}, and \eqref{lem-D-c4} after noting \eqref{lem-D-cond1} and \eqref{lem-D-cond2}. 

\medskip
The proof is complete. 
 \end{proof}

\begin{proposition} \label{pro-FTD}  Let $T > 0$ and $c>0$.  Let $(t_n)$ be an increasing sequence that converges to $T$ with $t_0 = 0$ and let $(\lambda_n), (\lambda_{1, n}) \subset \mR_+$ be increasing sequences. Assume that $\lambda_{1, n} \ge (2 + c) \lambda_n$. Set $s_0 = 0$ and $s_n = \sum_{k=0}^{n-1} \lambda_k(t_{k+1} - t_k)$ for $n \ge 1$.  There exists a positive constant $\gamma$ such that, if for large $n$, 
\be \label{pro-FTD-lambda1}
(t_{n+1} - t_n) \lambda_n \ge \gamma \lambda_n^{1/3},  \quad \mbox{ and }  \quad 
\lambda_{n+1} (t_{n+2} - t_{n+1}) \le (1 + 1/\gamma) \lambda_n (t_{n+1} - t_n), 
\ee
and 
\be \label{pro-FTD-lambda2}
\lim_{n \to + \infty} \frac{s_n}{n + \lambda_{n+1}^{1/3}} = + \infty, 
\ee
then there exists $\eps_0 > 0$ such that for all $y_0, \ty_0 \in L^2(0, L)$ with $\| y_0 \|_{L^2(0, L)}, \| \ty_0 \|_{L^2(0, L)} \le \eps_0$, there exists a unique pair $(y, \ty) \in X_T \times X_T$ such that, for $t_{n} \le t < t_{n+1}$ and $n \ge 0$, 
\be \label{cor-KdV-sys}
\left\{\begin{array}{c}
y_t + y_x + y_{xxx} + y y_x = 0 \quad \mbox{ in } (t_n, t_{n+1}), \\[6pt]
\ty_t + \ty_x + \ty_{xxx} + + 2 \lambda_n \ty - 
\lambda_{1, n} Q_n^{-1}(y - Q_n \ty) +   \ty y_x = 0 \mbox{ in } (t_n, t_{n+1}), \\[6pt]
y(\cdot, 0) = y (\cdot, L) = 0, y_x(\cdot, L) - y_x(\cdot, 0) = - \ty_x (\cdot, L) \mbox{ in } (t_n, t_{n+1}), \\[6pt]
\ty(\cdot, 0) = \ty (\cdot, L) = 0, \ty_x(\cdot, L) - \ty_x(\cdot, 0) = 0 \mbox{ in } (t_n, t_{n+1}), \\[6pt]
\end{array} 
\right. 
\ee
and
\be
y(0, \cdot) = y_0, \quad \ty(0, \cdot) = \ty_0 \mbox{ in } (0, L), 
\ee
where $Q_n = Q(\lambda_n)$ defined by \eqref{def-Q-KdV-PDE} with $\lambda = \lambda_n$. Moreover, we have, for $t_{n-1} \le t \le t_{n}$ and for $n \ge 1$,  
$$
\|y(t, \cdot) \|_{L^2(0, L)} \le e^{ - s_{n-1} + C n} \| y_0\|_{L^2(0, L)}. 
$$
for some positive constant $C$ independent of $n$, in particular, 
$$
y(t, \cdot) \to 0 \mbox{ in } L^2(0, L) \mbox{ as } t \to T_{-}. 
$$
\end{proposition}

\begin{remark} \rm There are sequences $(t_n)$, $(\lambda_n)$, $(\lambda_{1, n})$ which satisfy the conditions given in the above proposition, for example,  $t_n = T - T/n^2$, $\lambda_n = n^8$, and $\lambda_{1, n} = 2 \lambda_n$ for large $n$. 
\end{remark}

\begin{proof} Applying \Cref{pro-KdV-D} and  \Cref{lem-Qlambda-U,lem-Qlambda}, we have 
\be \label{pro-FTD-p1}
\| y(t_n, \cdot) \|_{L^2(0, L)} \le  e^{- 2 \lambda_{n-1} (t_n - t_{n-1}) + C (1+ \lambda_{n-1}^{1/3})} \| y(t_{n-1}, \cdot) \|_{L^2(0, L)} \mbox{ for } n \ge 1.
\ee
It follows that 
\be
\| y(t_{n}, \cdot) \|_{L^2(0, L)} \le e^{- s_{n-1} + C n} \| y_0 \|_{L^2(0, L)} \mbox{ for } n \ge 1. 
\ee
The conclusion now follows from \Cref{pro-KdV}. 

It is worth noting that \eqref{pro-FTD-lambda1} ensure the existence of $(y, \ty)$ by applying \Cref{pro-KdV-D} in the time interval $(t_n, t_{n+1})$ for all $n \ge 1$ (after a translation of time) since the condition \eqref{lem-KdVD-small} corresponding to the time interval $(t_n, t_{n+1})$, i.e., 
$$
\lambda_n^{2} e^{\frac{4}{3}\lambda_n (t_{n+1} - t_n)}) \| \|Q_n^{-1}\|^{2} \Big( \|y (t_n, \cdot) \|_{L^2(0, L)} + \| \ty (t_n, \cdot) \|_{L^2(0, L)} \Big) \le \eps,  
$$
is ensured if the following condition holds, for large $n$,  
$$
e^{\frac{4}{3} \lambda_n (t_{n+1} - t_n) + C \lambda_n^{1/3}} \| y(t_n, \cdot)\|_{L^2(0, L)} \le \eps. 
$$
if $\eps_0$ is sufficiently small. This holds by \eqref{pro-FTD-lambda1} and \eqref{pro-FTD-p1} 
\end{proof}



\bibliographystyle{amsref}




\end{document}